\documentclass[a4paper,12pt]{article}
\usepackage{amssymb}
\usepackage{amsmath}
\usepackage{array}
\usepackage{theorem}

\theoremstyle{change}
\theorembodyfont{\itshape}

\newcommand{\Pf}{{\em Proof}. }
\newcommand{\EPf}{\hfill$\square$}

\newcommand{\Span}[1]{\mbox{$<\!\!#1\!\!>$}}

\newcommand{\SU}[1]{\mbox{$\mathbf{SU}(#1)$}}
\newcommand{\U}[1]{\mbox{$\mathbf{U}(#1)$}}
\newcommand{\SP}[1]{\mbox{$\mathbf{Sp}(#1)$}}
\newcommand{\SO}[1]{\mbox{$\mathbf{SO}(#1)$}}

\newcommand{\Spin}[1]{\mbox{$\mathbf{Spin}(#1)$}}

\newcommand{\Lg}{\mbox{$\mathfrak g$}}
\newcommand{\Lk}{\mbox{$\mathfrak k$}}
\newcommand{\Lp}{\mbox{$\mathfrak p$}}
\newcommand{\La}{\mbox{$\mathfrak a$}}
\newcommand{\Lm}{\mbox{$\mathfrak m$}}
\newcommand{\Ln}{\mbox{$\mathfrak n$}}
\newcommand{\Lb}{\mbox{$\mathfrak b$}}

\newcommand{\Lu}{\mbox{$\mathfrak u$}}

\newtheorem{thm}{Theorem}

\newtheorem{prop}[thm]{Proposition}
\newtheorem{lem}[thm]{Lemma}

\begingroup\makeatletter\ifx\SetFigFont\undefined
\def\x#1#2#3#4#5#6#7\relax{\def\x{#1#2#3#4#5#6}}%
\expandafter\x\fmtname xxxxxx\relax \def\y{splain}%
\ifx\x\y   % LaTeX or SliTeX?
\gdef\SetFigFont#1#2#3{%
  \ifnum #1<17\tiny\else \ifnum #1<20\small\else
  \ifnum #1<24\normalsize\else \ifnum #1<29\large\else
  \ifnum #1<34\Large\else \ifnum #1<41\LARGE\else
     \huge\fi\fi\fi\fi\fi\fi
  \csname #3\endcsname}%
\else
\gdef\SetFigFont#1#2#3{\begingroup
  \count@#1\relax \ifnum 25<\count@\count@25\fi
  \def\x{\endgroup\@setsize\SetFigFont{#2pt}}%
  \expandafter\x
    \csname \romannumeral\the\count@ pt\expandafter\endcsname
    \csname @\romannumeral\the\count@ pt\endcsname
  \csname #3\endcsname}%
\fi
\fi\endgroup

\newcommand{\Aii}[2]{
\begin{picture}(1552,895)(1025,-4136)
\thinlines
\put(1201,-3961){\circle{336}}
\put(2401,-3961){\circle{336}}
\put(1369,-3961){\line( 1, 0){864}}
\put(1141,-3436){\makebox(0,0)[lb]{\smash{\SetFigFont{6}{7.2}{rm}#1}}}
\put(2341,-3436){\makebox(0,0)[lb]{\smash{\SetFigFont{6}{7.2}{rm}#2}}}
\end{picture}
}

\begingroup\makeatletter\ifx\SetFigFont\undefined
\def\x#1#2#3#4#5#6#7\relax{\def\x{#1#2#3#4#5#6}}%
\expandafter\x\fmtname xxxxxx\relax \def\y{splain}%
\ifx\x\y   % LaTeX or SliTeX?
\gdef\SetFigFont#1#2#3{%
  \ifnum #1<17\tiny\else \ifnum #1<20\small\else
  \ifnum #1<24\normalsize\else \ifnum #1<29\large\else
  \ifnum #1<34\Large\else \ifnum #1<41\LARGE\else
     \huge\fi\fi\fi\fi\fi\fi
  \csname #3\endcsname}%
\else
\gdef\SetFigFont#1#2#3{\begingroup
  \count@#1\relax \ifnum 25<\count@\count@25\fi
  \def\x{\endgroup\@setsize\SetFigFont{#2pt}}%
  \expandafter\x
    \csname \romannumeral\the\count@ pt\expandafter\endcsname
    \csname @\romannumeral\the\count@ pt\endcsname
  \csname #3\endcsname}%
\fi
\fi\endgroup

\newcommand{\Ai}[1]{
\begin{picture}(352,895)(1025,-4136)
\thinlines
\put(1201,-3961){\circle{336}}
\put(1141,-3436){\makebox(0,0)[lb]{\smash{\SetFigFont{6}{7.2}{rm}\mbox{$#1$}}}}
\end{picture}
}

\begin{document}

\title{Cycles of Bott-Samelson type for taut representations}
\author{Claudio Gorodski\footnote{\emph{Alexander von Humboldt Research
      Fellow} at the University of Cologne during the completion of this
      work.}\ \
      and Gudlaugur Thorbergsson}

\maketitle

\section{Introduction}

Ehresmann (\cite{Ehresmann}) introduced the 
Schubert varieties of the complex Grassmannians and showed 
that they give rise to a cell decomposition. Of the many 
generalizations of his result,
the most relevant for our work is the decomposition 
of the generalized real flag manifolds $G/P$ into Bruhat 
cells (see~\cite{Bruhat,Chevalley,HC,Warner}).
Here $G$ is a noncompact real semisimple 
Lie group with finite center and $P$ is a parabolic 
subgroup.

In their study of the Morse theory of symmetric spaces,
Bott and Samelson~(\cite{B-S}) came up with concrete cycles 
in the orbits of their isotropy representations which
represent a basis in $\mathbf Z_2$-homology.
It turns out that those orbits coincide with the
generalized real flag manifolds and the images of the 
cycles of Bott and Samelson equal the closures of the Bruhat 
cells (see~\cite{Hansen} and the appendix to this paper).
As an application of the cycles of Bott and Samelson, 
distance functions to orbits of isotropy representations
of symmetric spaces are perfect in the sense that the Morse equalities
are in fact equalities. We say that submanifolds with this property are 
taut and call representations taut if all of their orbits are taut. 
Because of their tautness, one sees that the Bruhat decomposition 
of the spaces $G/P$ is minimal in the sense that the number of cells 
in dimension $k$ is equal to the $k$th $\mathbf Z_2$-Betti number of
$G/P$. Notice that in the case $G$ is a complex group, the cells are all
even dimensional which makes the minimality of the Bruhat decomposition
trivial, but such easy arguments do not apply otherwise. 
Another Morse theoretic interpretation of the Bruhat cells
is their appearance as the unstable manifolds of suitable height functions 
on $G/P$ seen as orbits of isotropy representations of symmetric spaces,
see~\cite{Atiyah} and~\cite{Kocherlakota}. 

A construction similar to Bott and Samelson's was used by the second author 
in~\cite{Th2} to prove the tautness of complete proper Dupin
hypersurfaces. Hsiang, Palais and Terng were later able to
construct such cycles in isoparametric submanifolds~(\cite{HPT}), 
and Terng generalized that to weakly isoparametric submanifolds
in~\cite{Terng2} and to isoparametric submanifolds in Hilbert spaces
in~\cite{Terng3}. A further application of such constructions can be found
in~\cite{Th3} where a necessary topological condition is derived in order
for a manifold to admit a taut embedding. In~\cite{B-S} as well as
in the other papers, the cycles representing the homology of the submanifold 
(or orbit) are maps of iterated bundles 
of curvature surfaces (in general spheres) 
associated to the focal points along geodesics 
normal to that submanifold. Roughly speaking, as one moves along the normal
geodesic towards the submanifold, each focal point 
accounts for an iteration of the bundle, see Section~\ref{sec:BS}. 

We were able to show in~\cite{G-Th} that a taut irreducible 
representation of a compact connected Lie group is either 
the isotropy representation of a 
symmetric space or it is one of the following orthogonal representations
($n\geq2$): the $\mbox{(standard)}\otimes_{\mathbf R}\mbox{(spin)}$ 
representation of $\SO2\times\Spin9$; or the 
$\mbox{(standard)}\otimes_{\mathbf C}\mbox{(standard)}$
representation of
$\U2\times\SP n$; or the 
$\mbox{(standard)}^3\otimes_{\mathbf H}\mbox{(standard)}$
representation of $\SU2\times\SP n$. 
In this paper we will show how to adapt the construction of the cycles of Bott 
and Samelson to the orbits of these three representations.
As a result, they also admit explicit 
cycles representing a basis of their $\mathbf Z_2$-homology and, 
in particular, this provides another proof of their tautness. 

The main new technical difficulty that we encounter in our work is 
that focal points in the direction of normal vector fields parallel 
along some curvature circles
do not have constant multiplicity, 
which makes a modification in the construction of the cycles 
necessary to prevent the bundles from having some degenerate fibers. 
The modifications required by this 
``collapsing of focal points'' are only needed in finitely many points 
in one step of the construction and are achieved by a 
``cut and paste'' procedure, see
Section~\ref{sec:construction}. This strongly relies on a 
detailed knowledge of the geometry  
of the orbits of our representations, which is assembled
in Section~\ref{sec:reduction}. 

The first author wishes to thank the \emph{Alexander von Humboldt
Foundation} for its generous support and constant assistance
during the completion of this work.

\section{The method of Bott and Samelson}\label{sec:BS}
\setcounter{thm}{0}

Recall that a compact submanifold $M$ of an Euclidean space $V$ 
is called \emph{taut} with respect to the field of coefficients~$F$
if every (squared) distance function 
of a point $q\in V$, namely 
$L_q:M\to\mathbf R$ given by $L_q(x)=||q-x||^2$, which is a 
Morse function, has the minimum number of critical points 
allowed by the Morse inequalities with respect to~$F$. 
Throughout this paper we shall be assuming $F=\mathbf Z_2$ 
and dropping the reference to the field. 

Now assume that $M$ is a $G$-orbit of 
an orthogonal representation $\rho:G\to\mathbf O(V)$
of the compact connected Lie group $G$ and fix a $G$-regular point
$q\in V$ such that 
$L_q$ is a Morse function
and therefore has only finitely many 
critical points on $M$ with pairwise distinct critical values. 
For each critical point $p\in M$, 
Bott and Samelson constructed in~\cite{B-S} 
a compact manifold $\Gamma_p$ of dimension 
less than or equal to the index of $L_q$ at $p$ and a 
smooth map $h_p:\Gamma_p\to M$.
Under the assumption of variational completeness (see~\cite{B-S}
or~\cite{G-Th} for the definition of this concept)
they showed that the dimension of $\Gamma_p$ 
is equal to the index of $L_q$ at~$p$ and that  
\[ \oplus_p {h_p}_*:\bigoplus_p H_\lambda(\Gamma_p)\to H_\lambda(M) \] 
is an isomorphism, where $p$ runs through all the critical points of $L_q$ 
on $M$ such that the index of $L_q$ at~$p$ equals $\lambda$. 
This implies that the Morse inequalities for $L_q$ are equalities, i.~e.~$M$ is
taut. 

In fact, for a real number $c$ set 
\[ M^{c-}=\{x\in M:L_q(x)<c\}\quad\mbox{and}\quad
M^c=\{x\in M:L_q(x)\leq c\}. \]
We know that $M^c$ has the same homotopy type as $M^{c-}$ unless
$c$ is a critical value of $L_q$. Assume this is the 
case and let $p$ be the corresponding critical point, $L_q(p)=c$. 
Then $M^c$ has the homotopy type of $M^{c-}$ with a 
$\lambda$-cell $e_\lambda$ attached, where $\lambda$ is the 
index of $L_q$ at $p$. Consider the homology 
sequence of the pair $(M^c,M^{c-})$:
\begin{eqnarray*}
\lefteqn{\ldots\to\underbrace{H_{\lambda+1}(M^c,M^{c-})}_{\displaystyle=0}
\to H_\lambda(M^{c-})\to
  H_\lambda(M^c)} \\
\qquad&\to&\underbrace{H_\lambda(M^c,M^{c-})}_{\displaystyle=\mathbf Z_2}
\overset{\partial_*}{\to} H_{\lambda-1}(M^{c-})\to H_{\lambda-1}(M^c)
\to\underbrace{H_{\lambda-1}(M^c,M^{c-})}_{\displaystyle=0}\to\ldots
\end{eqnarray*}
By passing from $M^{c-}$ to $M^c$ the only possible changes 
in homology occur in dimensions $\lambda-1$ and $\lambda$. In the first case, 
the boundary $\partial e_\lambda$ of the attaching cell is a 
$(\lambda-1)$-sphere in $M^{c-}$ that does not bound a chain in $M^{c-}$,
so $e_\lambda$ has as boundary the nontrivial cycle $\partial e_\lambda$
in $M^{c-}$ and the map $\partial_*$ is not zero. 
In the second case, $\partial e_\lambda$ does bound a chain in~$M^{c-}$, 
which we cap with $e_\lambda$ to create a new nontrivial homology class
in~$M^c$, so $\partial_*$ is zero and $H_\lambda(M^c)\cong
H_\lambda(M^{c-})\oplus\mathbf Z_2$. Thus the Morse inequalities are
equalities when 
\setcounter{equation}{\value{thm}}
\stepcounter{thm}
\begin{equation}\label{eqn:surjectivity}
 H_{\lambda}(M^{c})\to H_{\lambda}(M^c,M^{c-1}) 
\end{equation}
is surjective, where $c$ is an arbitrary critical value of 
$L_q$ and $\lambda$ is the index of $L_q$ at the corresponding critical
point.
 
We proceed to describe the construction of Bott and Samelson
that yields the surjectivity of~(\ref{eqn:surjectivity}) under
the assumption of variational completeness. 
So fix $c$, $p$, $\lambda$ as above, let
$f_1,\ldots,f_r$ be the focal points of $M$ on the segment 
$\overline{qp}$ in focal distance decreasing order 
and let $m_1,\ldots,m_r$ be their respective
multiplicities. Note that $H=G_q$ is a principal isotropy group. 
We have that $\overline{qp}$ is perpendicular to $M$ at~$p$, therefore it
is also perpendicular to $Gq$ at~$q$, so that $H$ fixes the segment 
$\overline{qp}$ pointwise. 
Next let the $r$-fold product $H^r$ act on the product manifold 
$G_{f_1}\times\ldots\times G_{f_r}$ by the rule 
\[ \mathbf{g}\cdot
\mathbf{h}=(g_1h_1,\,h_1^{-1}g_2h_2,\,h_2^{-1}g_3h_3,\,\ldots,\,h_{r-1}^{-1}g_rh_r), \]
where $\mathbf{g}=(g_1,\ldots,g_r)\in G_{f_1}\times\ldots\times G_{f_r}$
and $\mathbf{h}=(h_1,\ldots,h_r)\in H^r$. Let $\Gamma_p$ be the quotient 
manifold under this action, namely, 
$\Gamma_p=G_{f_1}\times_H G_{f_2}\times_H\ldots\times_H G_{f_r}/H$,
and define $h_p:\Gamma_p\to M$ 
by $h[(g_1,\ldots,g_r)]=g_1\ldots g_rp$. 
It is immediate that $h_p$ is well-defined. 
Notice that $\Gamma_p$ is the total space of an iterated fiber bundle
that can be identified with the space of polygonal paths 
from~$p$ to $g_1\ldots g_rp$ with vertices~$g_rp$,
$g_{r-1}g_rp,\,\ldots,\,g_2\ldots g_{r-1}g_rp$ for 
$(g_1,\ldots,g_r)\in G_{f_1}\times\ldots\times G_{f_r}$. 
We compute:
\begin{eqnarray*}
 \dim\Gamma_p & = & (\dim G_{f_1}-\dim H)+\cdots+(\dim G_{f_r}-\dim H) \\
              & = & m_1 + \cdots + m_r \qquad\mbox{(by variational
              completeness, see Lemma~3.2 in~\cite{G-Th})} \\
              & = & \lambda \qquad\mbox{(by the index theorem of Morse).}
\end{eqnarray*}
Since $\Gamma_p$ is a compact manifold of dimension~$\lambda$, 
it follows that $H_\lambda(\Gamma_p)=\mathbf Z_2$.
Moreover, it is easy to see that $p\in h_p(\Gamma_p)\subset M^{c-}\cup\{p\}$
and that $h_p$ is an immersion near 
$(1,\ldots,1)\in G_{f_1}\times\ldots\times G_{f_r}$.
Now, locally in a Morse chart centered at~ $p$, the image
$h_p(\Gamma_p)$ is transversal to the ascending cell so that we can deform
it into the descending cell $e_\lambda$. 
Therefore 
${h_p}_*:H_\lambda(\Gamma_p)\to H_\lambda(M^c,M^{c-})$
is surjective. Finally, factorize
\[ H_\lambda(\Gamma_p)\to H_\lambda(M^c)
\to H_\lambda(M^c,M^{c-}) \]
to get the surjectivity of~(\ref{eqn:surjectivity}).

Let now $M$ be an arbitrary submanifold of an Euclidean space. 
A \emph{curvature surface} of $M$ is a submanifold $N$ 
such that there exists a parallel normal vector field 
$\xi$ of $M$ along $N$ such that 
the tangent space $T_r N$ is a full eigenspace
of the Weingarten operator $A_{\xi(r)}$ for all
$r\in N$. 

Turning back to the case of the $G$-orbit~$M$,
it follows easily from variational completeness
that for each focal point $f_i$ we have that 
$G_{f_i}p\subset M$ is a curvature surface of~$M$, 
so $\Gamma_p$ can be thought of
as an iterated bundle of
curvature surfaces. Those curvature surfaces will be spheres
if $q$ belongs to an open and dense subset of $V$. 

The cycle constructions in~\cite{HPT,Terng2,Terng3,Th2,Th3}
are all based on the observation that one does not need the action of the
group $G$ to define the iterated fiber bundle $\Gamma_p$ but only the
curvature surfaces themselves, which then must be assumed to be mutually
diffeomorphic in each step to guarantee that one indeed gets a fiber bundle. 
This follows in all those papers since the multiplicities of the
eigenvalues of the Weingarten operator $A_\xi$ are constant for 
parallel normal vector fields $\xi$ of~$M$ along curvature surfaces. 
However, that assumption is false for the three exceptional taut irreducible 
representations as we will see at the end of 
the next section. 

\section{The reduction principle}\label{sec:reduction}
\setcounter{thm}{0}

Let $\rho:G\to \mathbf O(V)$ be an orthogonal representation of a compact
Lie group $G$ which is not assumed to be connected. Denote by $H$ a fixed
principal isotropy subgroup of the $G$-action on $V$ and let $V^H$ be the
subspace of $V$ that is left pointwise fixed by the action of
$H$. Let $N$ be the normalizer of $H$ in $G$.
Then the group $N/H$ acts on $V^H$ with trivial principal
isotropy subgroup. Moreover, the following result is known
(\cite{G-S,Luna,LR,Schwartz,SS,Straume1}):

\begin{thm}[Luna-Richardson]\label{thm:LR}
The inclusion $V^H\to V$ induces a stratification preserving
homeomorphism between orbit spaces 
\setcounter{equation}{\value{thm}}
\stepcounter{thm}
\begin{equation}\label{eqn:homeo}
V^H/N\to V/G.
\end{equation}
\end{thm}

The injectivity of the map~(\ref{eqn:homeo}) means that 
$Np=Gp\cap V^H$ for $p\in V^H$. In particular, the $H$-fixed point 
set of a $G$-orbit is a smooth manifold.

Observe also that for a regular point $p\in V^H$
the normal space to the principal orbit $M=Gp$ at~$p$
is contained in $V^H$, because the slice representation at $p$
is trivial. More generally, we have:

\begin{lem}[\cite{G-Th}, Lemma~3.17]\label{lem:critical}
Let $p,q\in V^H$ and suppose that $q$ is a 
regular point for $G$. Consider $M=Gp$ and let $L_q:M\to\mathbf
R$ be the distance function. Then the critical set of $L_q$ is contained
in $M^H=M\cap V^H$, namely, $\mbox{Crit($L_q$)}=\mbox{Crit($L_q|_{M^H}$)}$.
\end{lem}

For $p\in V$ define $\nu_p=T_p(Gp)^\perp$ and for $p\in V^H$ 
define $\sigma_p=T_p(Np)^\perp\cap V^H$.
As a consequence of Lemma~\ref{lem:critical} we immediately get:

\begin{lem}
For $p\in V^H$, we have that $\sigma_p=\nu_p\cap V^H$.
\end{lem}

We also have: 

\begin{lem}[\cite{G-Th}, Lemma~3.16]\label{lem:weingarten}
For $p\in V^H$, $n\in\sigma_p$ and $M=Gp$ we have that the 
Weingarten operators satisfy $A_n^M|_{T_pM^H}=A_n^{M^H}$.
\end{lem}

Next fix $p\in V^H$ and consider the 
isotropy subgroup $N_p$. Then $N_p$ acts linearly on $\sigma_p$.
It is easy to apply the reduction principle as in Theorem~\ref{thm:LR}
to see that:

\begin{lem}\label{lem:slice}
The inclusion $\sigma_p\to\nu_p$ induces a 
stratification preserving homeomorphism between orbit spaces 
\[ \sigma_p/N_p\to\nu_p/G_p. \]
\end{lem}

In the following we specialize to the case where 
$\rho:G\to\mathbf O(V)$ is one 
of the three nonpolar taut irreducible 
representations, namely given by the 
following table ($n\geq2$; the last column will be explained later in this 
section):
\setlength{\extrarowheight}{0.3cm}
\[ \begin{array}{|c|c|c|c|c|}
\hline
case & G & \rho & \dim V=m & \mbox{$D$-$type$}\\
\hline
1&\SO2\times\Spin9 & \mbox{(standard)}\otimes_{\mathbf R}\mbox{(spin)} &
32& \Ai1\quad
\Ai6\quad \Ai7  \\
2&\U2\times\SP n & \mbox{(standard)}\otimes_{\mathbf C}\mbox{(standard)} &
8n&\quad\Ai1\quad \Ai2 \quad\Ai{\!\!\!\!\!\!\!4n-5}\quad\mbox{} \\
3&\SU2\times\SP n & \mbox{(standard)}^3\otimes_{\mathbf H}\mbox{(standard)} &
8n&\quad\Ai{\!\!\!\!\!\!\!4n-5} \quad\Aii11\quad\mbox{}\\
\hline
\end{array}\]

Then the cohomogeneity of $\rho$ is $3$, 
the dimension of $V^H$ is $4$ and $N/H$ consists of finitely many 
connected components each of which is a circle.
In particular the identity component $(N/H)^0$ is a circle 
group. Set $T=(N/H)^0$.
Then $N/H$ is a semidirect product $D\ltimes T$, where $D$ is a 
discrete group.   

Since $G$ acts linearly on $V$, it is enough to consider
its action on the unit sphere $S^{m-1}\subset V$.
Henceforth we view the orbits as submanifolds of $S^{m-1}$. 
The unit sphere $S^3\subset V^H$ is totally geodesically embedded
in $S^{m-1}$ as the the fixed point set $(S^{m-1})^H$. 

The reduced representation of $T\cong\SO2$ on 
$V^H\cong\mathbf R^4$ restricts to the Hopf action of 
$T$ on the unit sphere~$S^3\subset V^H$. 
For each $p\in S^3$, we shall call the great circle 
$Tp\subset Gp$ a \emph{special circle}. More generally, we shall 
call any $G$-translate of a special circle also a 
special circle. Notice that through each point 
$p\in S^{m-1}$ there exists a unique special circle, and 
that $S^3$ intersects every $G$-orbit in a finite number of 
circles, each of which is a special circle. 

An important consequence of the reduction principle for us 
is the technique which we call ``reduction of focal data'',
namely all the focal information about the $G$-orbits can be read off
the geometry of $S^3$. In fact, let $q\in S^{m-1}$ be a focal point 
of an orbit $Gp$. We can conjugate by an element of $G$ and then assume
that $q$ is a focal point of $Gp$ relative to~$p$ and that $p\in S^3$.
Moreover, it follows from Lemma~\ref{lem:slice} that

\begin{lem}
Let $q\in S^{m-1}$ be a focal point of $Gp$ relative to~$p\in S^3$. 
Then $q$ is $G_p$-conjugate to an element in $S^3$.
In particular, if $p$ is a regular point then $q\in S^3$.
\end{lem}

Next we discuss the multiplicities of the focal points of the orbits 
and distinguish between three different types of focal points.
It is useful to introduce the circle bundle 
\setcounter{equation}{\value{thm}}
\stepcounter{thm}
\begin{equation}\label{eqn:fibration}
\begin{array}{ccc}
T\approx S^1 & \rightarrow & S^3\subset V^H \\
 & &\!\!\!\!\!\!\!\!\!\!\!\!\!\!\!\!\!\!\!\!\!\!\!\eta\downarrow \\
 & & S^2\approx \mathbf CP^1
\end{array}
\end{equation}

We equip the base space $S^2$ with the quotient metric. Then
$S^2$ is a metric sphere of radius $1/2$ and 
the discrete group~$D$ acts there by isometries.
The $D$-singular set is a union of great circles,
and each $D$-singular great circle $C$ comes with a multiplicity, namely
the difference between the dimension of a principal $G$-orbit and the
dimension of the orbit $Gq$ where $q$ is a point 
in~$\eta^{-1}(C)$ not contained in the preimage of any other
$D$-singular great circle. 
We write the multiplicity of $C$ next to the vertex 
of the diagram of $D$ which corresponds to $C$ (or to a $D$-conjugate of
$C$), see the last column in the table. 

For fixed $p\in S^3$ and $n\in\sigma_p\cap T_pS^3$ a unit vector, let 
$\gamma_t(s)=(\cos s)tp + (\sin s)n$, where $s\in\mathbf R$, $t\in T$,
be the normal geodesic to $Gp$ in the direction of $n\in\sigma_{tp}$.
Notice that $\gamma_t$ is an horizontal curve with respect 
to the Riemannian submersion~$\eta$ and that $n$ as a point 
in~$S^3$ satisfies $n\in Gp$ (because $\eta(n)=-\eta(p)$ is 
$D$-conjugate to $\eta(p)$). 

If $p\in S^{m-1}$ and $q$ is a focal point of $Gp$ relative to~$p$
of multiplicity $k>0$, we shall call $q$ a focal point of 
\emph{standard type} (resp.~\emph{mixed type}, \emph{special type}) if 
$k=\dim G_q p$ (resp.~$k>\dim G_q p>0$, $\dim G_q p=0$). 
Any focal point falls into one of these three cases. 
Notice that the first case occurs precisely 
when the focal point satisfies the 
condition in the definition of variational completeness.

\begin{prop}\label{prop:focal}
With the above notation:
\begin{enumerate} 
\item The orbit spaces $S^{m-1}/G\cong S^2/D$. 
\item If $p\in S^3$ 
is a regular point, then the focal points of standard type of
$Gp$ along the normal geodesic segment $\gamma_1|_{[0,\pi)}$ are precisely the 
points $\gamma_1(s)$, $0<s<\pi$, such that $\eta(\gamma_1(s))$ 
belongs to a $D$-singular great circle $C$ in~$S^2$. The 
multiplicity of a standard focal point $\gamma_1(s)$
is the sum of the multiplicities of the 
$D$-singular great circles which pass through $\eta(\gamma_1(s))$.
There is precisely one focal point of special type along 
$\gamma_1|_{[0,\pi)}$, namely $n\in S^3$, its multiplicity is one, 
its focal distance is $\pi/2$ and the associated line of curvature is 
$Tp$. There are no focal points of mixed type.
\item If $p\in S^3$ is a singular point and 
$d\eta_p(n)$ is not tangent to a $D$-singular great circle in~$S^2$, 
then the focal points of standard type of
$Gp$ along the normal geodesic segment $\gamma_1|_{[0,\pi)}$ are precisely the 
points $\gamma_1(s)$, $0<s<\pi$ and $s\neq\pi/2$, 
such that $\eta(\gamma_1(s))$ 
belongs to a $D$-singular great circle $C$ in~$S^2$. The 
multiplicity of a standard focal point $\gamma_1(s)$
is the sum of the multiplicities of the 
$D$-singular great circles which pass through $\eta(\gamma_1(s))$.
There is precisely one focal point of mixed type along along 
$\gamma_1|_{[0,\pi)}$, namely $n\in S^3$, its multiplicity is 
$\dim G_n p +1$ and its focal distance is $\pi/2$. 
There are no focal points of special type.
\end{enumerate}
\end{prop}

\Pf Theorem~\ref{thm:LR} and~(\ref{eqn:fibration}) 
imply that $S^{m-1}/G\cong S^3/N\cong S^2/D$. This gives~(a).
If $p$ is a regular point, then each point $\gamma_1(s)$, $0<s<\pi$, 
such that $\eta(\gamma_1(s))$ 
belongs to a $D$-singular great circle $C$ in~$S^2$ is a focal point
of $Gp$ of multiplicity greater than or equal to the 
sum of multiplicities of the $D$-singular great circles which 
pass through $\eta(\gamma_1(s))$. It follows from a counting 
argument based on data from the table that the sum of the multiplicities 
of the focal points of $Gp$ along the normal geodesic
segment $\gamma_1|_{[0,\pi)}$ obtained in this way is 
greater than or equal to $\dim Gp-1$, none of which has
a focal distance equal to $\pi/2$. Now $n\in\sigma_{tp}$ 
for all $t\in T$, so $n\in S^3$ is a focal point of $Gp$
with focal distance $\pi/2$. Since the sum of the multiplicities 
of the focal points to $Gp$ along the geodesic segment $\gamma_1|_{[0,\pi)}$
must be equal to $\dim Gp$, we have that there are 
no other focal points and~(b) follows.
If $p$ is a singular point, we observe: $n$ is also a singular point;
an element in $G_n$ that fixes $p$ must also fix the geodesic 
segment between $n$ and $p$. Now 
the proof of~(c) is similar to the above.  
\EPf

\medskip

In the cycle constructions in~\cite{HPT,Terng2,Terng3,Th2,Th3}
and implicitly in~\cite{B-S}, 
it is essential that the multiplicities of the eigenvalues of $A_\xi$ be
constant for $\xi$ a parallel normal vector field along a curvature
surface; see the end of Section~\ref{sec:BS}. We will now see that we do 
not have such constancy for certain parallel normal vectors along the 
special circles. Notice that the nonconstancy of multiplicities of
eigenvalues of $A_\xi$ is equivalent to the nonconstancy of 
multiplicities of focal points in the direction of~$\xi$, which in turn is
equivalent to saying that there is collapsing of focal points in the 
direction of~$\xi$.  

In fact, the differential $d\eta_p:\sigma_p\cap T_pS^3\to T_{\eta(p)}S^2$ 
is an isometry and $\tau(t)=d\eta_{tp}\circ P_t\circ d\eta_p^{-1}$
for $t\in T$ defines the standard representation
$\tau:T\cong\SO2\to T_{\eta(p)}S^2\cong\mathbf R^2$,
where $P_t$ denotes parallel transport along the special circle 
from $p$ to~$tp$. Moreover, 
\[ \eta(\gamma_t(s))=(\cos s)\eta(p)
                               +(\sin s)\tau(t)d\eta_p(n),
\]
so the distribution of focal points of standard type on the geodesic
segment $\gamma_t|_{[0,\pi/2]}$ from $tp$ to $n$ 
depends on $t\in T$, but is well controlled
by the intersections of $\eta\circ\gamma_t|_{[0,\pi/2]}$ with 
the $D$-singular set in $S^2$. 
In particular, for a finite number of values $t\in T$ there is collapsing
of focal points in the direction of~$n$ viewed as a parallel 
normal vector field along the special circle~$Tp$. 

\medskip

The existence of focal points other than of standard type precludes 
variational completeness (see Lemma~3.1 in~\cite{G-Th}). Therefore:

\begin{thm}[\cite{G-Th}]
The representations listed in the table are not variationally complete.
\end{thm}

\section{The construction of the cycles}\label{sec:construction}
\setcounter{thm}{0}

According to the discussion in the previous section, the representations
listed in the table there
fail to be variationally complete just because
the Jacobi field $J$ along $\gamma_1$ defined by
$J(s)=(\cos s)v$, where $v$ is a vector tangent to the 
special circle $Tp$ at~$p$, 
is not the restriction of a Killing 
vector field induced by the representation. 
The ``variational cocompleteness''  
of the principal orbits of 
these representations in the sense of~\cite{Hsiang} is one, 
and in this regard they are very close to being 
variationally complete. In this section we want
to adapt the construction described 
in Section~\ref{sec:BS} to obtain cycles 
of Bott-Samelson type for the orbits of these 
nonvariationally complete representations. 

So let $M$ be a $G$-orbit of one of the representations listed in the table,
fix a $G$-regular point 
$q\in V$ such that the distance function $L_q:M\to\mathbf R$ is a Morse
function and let $p\in M$ be a critical point of $L_q$ 
such that $L_q(p)=c$ and the index of $L_q$ at~$p$ is $\lambda$. 
We want to construct a compact $\lambda$-dimensional manifold $\Gamma_p$ 
and a smooth map $h_p:\Gamma_p\to M$ such that 
$p\in h_p(\Gamma_p)\subset M^{c-}\cup\{p\}$ and $h_p$ is an immersion 
near $h_p^{-1}(p)$. As in Section~\ref{sec:BS}, this will
imply the surjectivity of~(\ref{eqn:surjectivity})
and therefore the tautness of~$M$.  

Before we begin, we introduce the 
following notation for future reference.
If $q_1,\ldots,q_u\in S^3$, we denote by
$W_{q_1,\ldots,q_u}$ the quotient manifold
$G_{q_1}\times_H G_{q_2}\times_H\ldots\times_H G_{q_u}/H$,
see Section~\ref{sec:BS}.

\textit{We start the construction with the case where $M$ is a principal
orbit.} If all the focal points to~$M$
along the normal geodesic segment $\overline{qp}$
are of standard type, we follow the same procedure
as in Section~\ref{sec:BS} in order to define $\Gamma_p$ and $h_p$. 

Consider next the case where the furthermost focal point, 
say $f_1$, is of special type. Then all the other focal 
points on the normal geodesic segment are of standard type. 
It follows from Proposition~\ref{prop:focal} and the discussion 
in Section~\ref{sec:reduction} that there is a finite subset 
$\{t_0,\ldots,t_{l-1}\}\subset T$ such that the multiplicities 
of the focal points to~$M$ along the geodesic segment 
$\overline{q(tp)}$ are locally constant for 
$t\in T\setminus\{t_0,\ldots,t_{l-1}\}$. Moreover, we may perturb 
$q$ if necessary in order to have $t_i\neq 1$ for $i=0,\ldots,l-1$.
Note that for all $t\in T\setminus\{t_0,\ldots,t_{l-1}\}$
we can specify the focal points along $\overline{q(tp)}$ as
$f_1,f_2(t),\ldots,f_r(t)$ in focal distance decreasing order. 
Next we assume the $t_i$'s are ordered so that 
they partition the circle $T$ in closed arcs 
$[t_0,t_1],[t_1,t_2],\ldots,[t_{l-2},t_{l-1}],[t_{l-1},t_0]$.
Fix $[t_{i-1},t_i]$ (indices~$i$ modulo~$l$) and
define a fiber bundle $\Gamma_p^i\to[t_{i-1},t_i]$ 
as follows. For $t\in(t_{i-1},t_i)$, the isotropy 
subgroups $G_{f_j(t)}$ for $j=2,\ldots,r$ are well-defined 
and we can take the fiber $\Gamma_p^i|_t$ of 
$\Gamma_p^i$ over $t\in(t_{i-1},t_i)$ to be 
$W_{f_2(t),\ldots,f_r(t)}$. 
We want to define $G_{f_j(t_{i-1}^+)}$ and $G_{f_j(t_i^-)}$
for $j=2,\ldots,r$. 
Note that as $t\to t_{i-1}^+$ (resp.~$t\to t_i^-$)
two or more of the focal points of $M$ 
relative to $tp$ collapse into a focal point relative 
to $t_{i-1}p$ (resp.~$t_ip$) 
whose multiplicity is equal to the sum of the
multiplicities of the collapsing focal points.
If $j$ is an index for which the focal point $f_j(t)$ is not collapsing
as $t\to t_{i-1}^+$ (resp.~$t\to t_i^-$) we simply take 
$G_{f_j(t_{i-1}^+)}$ (resp.~$G_{f_j(t_i^-)}$) to be the isotropy 
subgroup $G_{f_j(t_{i-1})}$ (resp.~$G_{f_j(t_i)}$).
Otherwise, suppose that $j$ is an index for which 
$f_j(t)$ collapses as $t\to t_{i-1}^+$ (resp.~$t\to t_i^-$) 
into a focal point $\bar f$. Then 
$\eta(f_j(t))$ belongs to a $D$-singular great circle $C$
which also contains $\eta(\bar f)$ (here $t\in(t_{i-1},t_i)$). 
Define the subgroup $G_{\tilde C}\subset G_{\bar f}$ to be the 
stabilizer of the $\eta$-horizontal lift of~$C$ through $\bar f$, 
which we call $\tilde C$.
Alternatively, $G_{\tilde C}$ could be defined as 
a certain isotropy subgroup relative to the slice representation 
at~$\bar f$. Since $t\mapsto f_j(t)$ and $\tilde C$ are both smooth lifts 
of $C$ for $t\in(t_{i-1},t_i)$, it follows that there exists a
smooth curve $n_j(t)\in T$ such that 
$G_{f_j(t)}=n_j(t)G_{\tilde C}n_j^{-1}(t)$ for $t\in(t_{i-1},t_i)$. 
Moreover, as $t\to t_i^+$ (resp.~$t\to t_{i+1}^-$) 
we have that 
$f_j(t)\to\bar f$, so $\lim_{t\to t_{i-1}^+}n_j(t)=1$ 
(resp.~$\lim_{t\to t_i^-}n_j(t)=1$)
and $n_j$ can 
be smoothly extended to $[t_{i-1},t_i)$
(resp.~$(t_{i-1},t_i]$). 
Let $G_{f_j(t_{i-1}^+)}=G_{\tilde C}$ 
(resp.~$G_{f_j(t_i^-)}=G_{\tilde C}$).
Now we can define
the fiber $\Gamma_p^i|_{t_{i-1}}$ (resp.~$\Gamma_p^i|_{t_i}$)
of $\Gamma_p^i$ over $t_{i-1}$ (resp.~$t_i$) to be 
$W_{f_2(t_{i-1}^+),\ldots,f_r(t_{i-1}^+)}$, 
(resp.~$W_{f_2(t_i^-),\ldots,f_r(t_i^-)}$).
The previous discussion shows that $\Gamma_p^i$ is a smooth 
fiber bundle over $[t_{i-1},t_i]$.
We also define a smooth map $h_p^i:\Gamma_p^i\to M$ by
$h_p^i[(g_2,\ldots,g_r)]=g_2\ldots g_rp$, where
$(g_2,\ldots,g_r)\in\Gamma_p^i|_t$. 

For each $i:1,\ldots,l$,
note that $h_p^i(\Gamma_p^i|_{t_i})=h_p^{i+1}(\Gamma_p^{i+1}|_{t_i})$
and denote this subset of $M$ by $A_i$. 
The last step in the 
construction is the definition of a `glueing' diffeomorphism
$\psi_i:\Gamma_p^i|_{t_i}\to\Gamma_p^{i+1}|_{t_i}$ 
and a `correcting' diffeomorphism 
$\phi_i:M\to M$ for each 
$i=1,\ldots,l$ such that $\phi_{i+1}(A_i)=A_i$ and 
\setcounter{equation}{\value{thm}}
\stepcounter{thm}
\begin{equation}\label{eqn:h}
h_p^i|_{\Gamma_p^i|_{t_i}}=\phi_{i+1}h_p^{i+1}\psi_i
\end{equation}
for $i=1,\ldots,l$ (indices~$i$ modulo~$l$).
We also require the `cycle condition'
\setcounter{equation}{\value{thm}}
\stepcounter{thm}
\begin{equation}\label{eqn:cycle}
\phi_1\phi_2\cdots\phi_l=1.
\end{equation}
Then we will be able to define $\Gamma_p$ as the resulting fiber bundle
over $S^1$ and $h_p:\Gamma_p\to M$ as the map that restricts to 
$\phi_2\phi_3\cdots\phi_i h_p^i$ on each closed subset $\Gamma_p^i$. 

According to the discussion in Section~\ref{sec:reduction},
for each $t_i$ we have either one collapsing of   
two focal points (\emph{double collapse};
occurs six times for each representation listed
in the table) or one collapsing of 
three focal points (\emph{triple collapse};
occurs twice for the representation listed in the table
under number~3).
We discuss these two possibilities separately.

\bigskip
\textsc{Double collapse}
\bigskip

Suppose that for some $j=2,\ldots,r-1$ we have that
$f_j(t)$ and $f_{j+1}(t)$ collapse together into a focal point $\bar f$
as $t\to t_i$. 
In this case we have that $G_{f_j(t_i^+)}=G_{f_{j+1}(t_i^-)}$
and $G_{f_{j+1}(t_i^+)}=G_{f_j(t_i^-)}$.
It is not difficult to see that 
$G_{f_j(t_i^+)}\cap G_{f_{j+1}(t_i^+)}=H$ 
and that $(G_{f_j(t_i^+)},G_{f_{j+1}(t_i^+)})$ is a factorization of
the isotropy subgroup $G_{\bar f}$, 
that is
$G_{\bar f}=G_{f_j(t_i^+)}\cdot G_{f_{j+1}(t_i^+)}$. 
It follows that the multiplication map 
$G_{f_j(t_i^+)}\times G_{f_{j+1}(t_i^+)}\to G_{\bar f}$
induces a diffeomorphism
\[ G_{f_j(t_i^+)}\times_H G_{f_{j+1}(t_i^+)}\approx G_{\bar f}.\]
Clearly there is a similar diffeomorphism
$ G_{f_j(t_i^-)}\times_H G_{f_{j+1}(t_i^-)}\approx G_{\bar f}$. 
Now the composed diffeomorphism 
$G_{f_j(t_i^-)}\times_H G_{f_{j+1}(t_i^-)}\to 
G_{f_j(t_i^+)}\times_H G_{f_{j+1}(t_i^+)}$
extends trivially to a diffeomorphism 
\[ \psi_i:W_{f_2(t_i),\ldots,f_j(t_i^-),f_{j+1}(t_i^-),
\ldots,f_r(t_i)}
\to W_{f_2(t_i),\ldots,f_j(t_i^+),f_{j+1}(t_i^+),
\ldots,f_r(t_i)} \] 
which satisfies~(\ref{eqn:h}) if we take $\phi_{i+1}=1$. 
Observe that this choice of $\phi_{i+1}$ leaves 
unaffected the cycle condition~(\ref{eqn:cycle}). 

\bigskip
\textsc{Triple collapse}
\bigskip

Suppose that for some $j=2,\ldots,r-2$ we have that
$f_j(t)$, $f_{j+1}(t)$ and $f_{j+2}(t)$ collapse together 
into a focal point $\bar f$ as $t\to t_i$. 
In this case we have that $G_{f_j(t_i^+)}=G_{f_{j+2}(t_i^-)}$,
$G_{f_{j+1}(t_i^+)}=G_{f_{j+1}(t_i^-)}$
and $G_{f_{j+2}(t_i^+)}=G_{f_j(t_i^-)}$.
This is the case of representation number~3 in the table,
so the group $D$ is of $A_1\times A_2$-type, namely
$D$ is isomorphic to $\mathbf Z_2\oplus\mathbf D_3$
where $\mathbf D_3$ is the dihedral group of degree~$3$.
Note that the $N$-isotropy subgroup at~$\bar f$ equals
$\mathbf D_3$, and let $w^\circ\in\mathbf D_3$ 
be the unique element of order $2$
which conjugates $G_{f_j(t_i^+)}$ into $G_{f_{j+2}(t_i^+)}$. 
Then $w^\circ$ normalizes $G_{f_{j+1}(t_i^+)}$ and $H$, so 
the diffeomorphism 
\[ G_{f_j(t_i^-)}\times_H G_{f_{j+1}(t_i^-)}\times_H G_{f_{j+2}(t_i^-)}\to
G_{f_j(t_i^+)}\times_H G_{f_{j+1}(t_i^+)}\times_H G_{f_{j+2}(t_i^+)} \]
given by 
\[ (g_1,g_2,g_3)\mapsto(w^\circ g_1 w^\circ,
w^\circ g_2 w^\circ,w^\circ g_3 w^\circ) \]
is well-defined and extends trivially to a diffeomorphism 
\[ \psi_i:W_{f_2(t_i),\ldots,f_j(t_i^-),f_{j+1}(t_i^-),f_{j+2}(t_i^-),
\ldots,f_r(t_i)}
\to W_{f_2(t_i),\ldots,f_j(t_i^+),f_{j+1}(t_i^+),f_{j+2}(t_i^+),
\ldots,f_r(t_i)} \] 
which satisfies~(\ref{eqn:h}) if we take $\phi_{i+1}$ to be 
conjugation by~$w^\circ$. Observe that the triple collapse 
occurs precisely twice for this representation, and it is 
easy to see that in both 
instances the element $w^\circ$ is the same. Since $w^\circ$
has order two, condition~(\ref{eqn:cycle}) is also satisfied. 

\medskip

The above discussion finishes the construction for the case where 
the furthermost focal point on the 
normal geodesic segment $\overline{qp}$
is of special type. We now consider the case 
where there is a focal point of special type on that segment which is not the
furthermost focal point, say it is $f_r$, 
the $r$-th focal point in focal distance decreasing order
for $r\geq2$. Let $f_{r-1}$ be the previous focal point
and choose an interior point in the segment $\overline{f_{r-1}f_r}$
to be denoted~$q'$.  
We have that the furthermost focal point to~$M$
along the normal geodesic segment~$\overline{q'p}$ is~$f_r$, hence 
of special type, so 
we may apply the previous case and obtain a 
cycle~$h_p^0:\Gamma_p^0\to M$ for this situation.  
Note that~$H$ acts naturally on the left on the manifold~$\Gamma_p^0$
(this is just left multiplication for representations number~1 and~2;
for representation number~3, note that
the two triple collapses correspond to points in the special circle
which cut it into two halves, and that the element
$h\in H$ acts by left multiplication by $h$ on the fibers lying over 
one half of the special circle and acts by 
left multiplication by $w^\circ hw^\circ$ on the fibers lying over
the other half).
The~$r-1$ first focal points~$f_1,\ldots,f_{r-1}$ along $\overline{qp}$ 
are all of standard type. We define~$\Gamma_p=G_{f_1}\times_H 
\ldots\times G_{f_{r-1}}\times_H \Gamma_p^0$ and~$h_p:\Gamma_p\to M$
by $h_p[(g_1,\ldots,g_{r-1},\mathbf g_0)]=g_1\ldots g_{r-1}h_p^0(\mathbf
g_0)$, where $(g_1,\ldots,g_{r-1})\in G_{f_1}\times\ldots\times G_{f_{r-1}}$,
$\mathbf g_0\in\Gamma_p^0$. 
This concludes the discussion for the case of a principal orbit~$M$. 

In the case where $M$ is not a principal orbit, 
the procedure is similar. We just remark that when the furthermost 
focal point along the normal geodesic segment $\overline{qp}$, 
say $f_1$, is of mixed type, then $f_1=n$ and one does ``cut and paste'' 
to construct $\Gamma_p$ as a fiber bundle over~$S^1$ 
using typical fibers of the form $W_{f_1(t),f_2(t),\ldots,f_s(t)}$. 
This completes the construction 
of the cycles of Bott-Samelson type
for the exceptional taut irreducible representations. 

\bigskip
\noindent\textbf{Questions } 1.~Are there orthogonal representations 
which have principal orbits with ``variational cocompleteness'' 
one (see~\cite{Hsiang}) without being taut?

\noindent 2.~Is it possible to generalize the construction 
given in this paper and find explicit cycles for the 
$\mathbf Z_2$-homology of an arbitrary taut submanifold~$M$ of 
an Euclidean space?

\noindent 3.~Assume the answer to the question in~(2) to be positive, 
say for some subclass of taut submanifolds at least. Is it then 
possible to give an upper bound for the number of cycles 
and hence estimate the $\mathbf Z_2$-Betti numbers of~$M$?
One might conjecture the upper bound for the $k$th Betti number
to be $\binom{k}{n}$ where $n$ is the dimension of~$M$.

\appendix

\section*{Appendix: The Bruhat cells versus the Bott-Samelson cycles}
\setcounter{thm}{0}
\setcounter{section}{1}
\renewcommand{\thesection}{\Alph{section}}

In this appendix we prove that the closures of the Bruhat cells coincide
with the images of the Bott-Samelson cycles in a real generalized flag
manifold $G/P$. We are not aware of any such proof in the literature, but the
special case where $G$ is a complex group and $P$ is a Borel subgroup
is treated by Hansen in~\cite{Hansen}. Our proof is analogous to his. 

We follow closely (but not strictly) 
the setting and notation from Sections~1.1
and~1.2 in~\cite{Warner}. 
Let $G$ be a noncompact real connected semisimple Lie group 
with finite center and denote by~$\Lg$ its Lie algebra. 
Consider a fixed Cartan decomposition $\Lg=\Lk+\Lp$
with corresponding Cartan involution~$\theta$ 
and equip~$\Lg$ with the real inner product~$(X,Y)_\theta=-B(X,\theta Y)$
where~$X,Y\in\Lg$ and~$B$ is the Killing form of~$\Lg$.  
Let $\La$ be a maximal Abelian subspace in~$\Lp$. 
Then we have the real orthogonal restricted root decomposition
$\Lg = \Lm + \La + \sum_{\lambda\in\Sigma} \Lg_\lambda$
where $\Lm$ is the centralizer of $\La$ in~$\Lk$. Introduce an 
order in the dual $\La^*$ and let $\Sigma^+$ (resp.~$\Sigma^-$)
denote the set of positive (resp.~negative)
restricted roots. Then $\Ln^+=\sum_{\lambda\in\Sigma^+} \Lg_\lambda$
and $\Ln^-=\theta(\Ln^+)=\sum_{\lambda\in\Sigma^-} \Lg_\lambda$
are nilpotent subalgebras of $\Lg$. Let $K$, $A$, $N^+$, $N^-$  
be respectively the analytic subgroups of $G$ with Lie algebras 
$\Lk$, $\La$, $\Ln^+$ and~$\Ln^-$. 
Then $K$, $A$, $N^+$ and $N^-$ are closed subgroups of $G$ and 
we have the Iwasawa decomposition $G=KAN^+$. Let $M$ be the
centralizer of $A$ in $K$ and let $M^*$ be the normalizer of $A$
in~$K$. Then the Weyl group of $(\Lg,\La)$ generated by reflections on the 
singular hyperplanes in~$\La$ (which are the
kernels of the reduced restricted roots) is identified 
with $M^*/M$ and $B=MAN^+$ is a Borel subgroup of $G$. 
Define also $U_w^+=N^+\cap wN^-w^{-1}$ where $w\in W$ (we abuse the
notation in the evident manner). Then we have 
(see~\cite{Bruhat,Chevalley,HC,Warner}):

\begin{thm}\label{thm:B-HC} 
The quotient manifold $G/B$ is the disjoint union 
of the $N^+$-orbits $N^+w^{-1}B$ for $w\in W$. Moreover, for each $w\in W$,
the map $U_{w^{-1}}^+\to N^+w^{-1}B$, $v\mapsto vw^{-1}B$ is a diffeomorphism.
\end{thm}

Theorem~\ref{thm:B-HC} admits the following generalization.
Let $\Upsilon\subset\Sigma^+$ be the system of simple roots.
The cardinality $l$ of $\Upsilon$ is equal to the dimension of $\La$
and $l$ is the split-rank of $\Lg$. Fix 
a subset $\Theta$ of $\Upsilon$ and denote by $\Span{\Theta}$ the set
comprised of those $\lambda$ in $\Sigma$ which are linear combinations
of the elements of $\Theta$; we also agree to write 
$\Span{\Theta}^{\pm}$ for $\Sigma^{\pm}\cap\Span{\Theta}$. 
Put $\Ln^{\pm}(\Theta)=\sum_{\lambda}\Lg_\lambda$, 
$\lambda\in\Span{\Theta}^{\pm}$
and $\Ln_{\Theta}^+=\sum_\lambda \Lg_\lambda$, 
$\lambda\in\Sigma^+\setminus\Span{\Theta}^+$; $\Ln_\Theta^-=\theta(\Ln_\Theta^+)$. 
Put also $\La(\Theta)=\sum_\lambda\mathbf RH_\lambda$, 
$\lambda\in\Span{\Theta}^+$, where $H_\lambda\in\La$ is the coroot
associated to $\lambda$, and let $\La_\Theta$ be the orthogonal complement 
of $\La(\Theta)$ in $\La$. 
Set $\Lp_\Theta=\Lb+\Ln^-(\Theta)$ where $\Lb=\Lm+\La+\Ln^+$ 
is the Lie algebra of $B$. Then the space $\Lp_\Theta$ is the normalizer
of $\Ln^+_\Theta$ in $\Lg$. Hence $\Lp_\Theta$ is a Lie subalgebra 
of $\Lg$ and the corresponding analytic subgroup $P_\Theta$ of $G$ is 
closed. We have that $P_\Theta$ is a parabolic subgroup of $G$, and 
the $2^l$ subgroups constructed in this way for each possible choice of 
$\Theta$ are the so-called standard parabolic subgroups of $G$ 
(with respect to our choice
of Iwasawa decomposition). 
It can be shown that $P_\Theta$
admits the Langlands decomposition $P_\Theta=M_\Theta(K)AN^+$ 
where $M_\Theta(K)$ is the $K$-centralizer of $\La_\Theta$.
For each $w\in W$ 
write $\Sigma_w^+=w\Sigma^-\cap\Sigma^+$ and 
denote by $W_u$ the set
of those $w\in W$ such that $\Sigma_w^+\cap\Span\Theta^+=\varnothing$. 
Now we have:

\begin{thm}\label{thm:B-K} 
The quotient manifold $G/P_\Theta$ is the disjoint union 
of the $N^+$-orbits $N^+w_u^{-1}P_\Theta$ for $w_u\in W_u$. 
Moreover, for each $w_u\in W_u$,
the map $U^+_{w_u^{-1}}\to N^+w_u^{-1}P_\Theta$, 
$v\mapsto vw_u^{-1}P_\Theta$ is a diffeomorphism.
\end{thm}

In~\cite{Warner}
Theorems~\ref{thm:B-HC} and~\ref{thm:B-K} are 
respectively accredited to Bruhat -- Harish-Chandra
and~Borel -- Kostant.  
Note that Theorem~\ref{thm:B-HC}
is the particular instance of Theorem~\ref{thm:B-K} 
corresponding to $\Theta=\varnothing$. 
Note also that Theorem~\ref{thm:B-K}
gives a cell decompositions for $G/P_\Theta$
which turns it into a finite CW-complex.

Next let $X$ be a Riemannian symmetric space with no Euclidean 
factors which is also assumed to be noncompact, since 
the isotropy representation of a 
compact symmetric space coincides with that of its noncompact dual.
Then we can find $G, K$ as above such that $X=G/K$,
namely $G$ is the connected component of the isometry group of
$X$ and $K$ is the isotropy subgroup at some point.  
Now the isotropy representation of $X$ at $1\cdot K$ is equivalent 
to the adjoint representation of $K$ on $\Lp$. 
Fix once and for all $p\in\La$. 
The isotropy subgroup $K_p=M_\Theta(K)$ where
$\Theta$ consists of those simple restricted roots 
which vanish at~$p$. Thus we can identify the $K$-orbit through~$p$
with $K/M_\Theta(K)=G/P_\Theta$. 

Let $q\in\La$ be a regular element and consider the distance function 
$L_q:Kp\to\mathbf R$. Then $L_q$ is a Morse function whose critical points 
are precisely the Weyl translates $wp$ for $w\in W$. Therefore 
the Bott-Samelson cycles associated to $L_q$ can be parametrized 
by $W/W_p$. For a given $w\in W$, 
consider the segment joining $q$ and $wp$ and suppose 
that this segment cuts across the singular hyperplanes corresponding to the 
reduced restricted roots $\lambda_1,\ldots,\lambda_k$ in this order.
Then it is clear that $w$ and $s_k\cdots s_1$ are $W_p$-conjugate
where $s_i$ is the 
reflection in the singular hyperplane $\{H\in\La:\lambda_i(H)=0\}$;
let $K_i$ denote the closed subgroup of $K$ which is the stabilizer 
of that singular hyperplane. Denote the coset $wW_p$ by $\bar w$. 
The Bott-Samelson cycle associated to $\bar w$ may be assumed to be
$h_{\bar w}:\Gamma_{\bar w}\to K/K_p$, where 
$\Gamma_{\bar w}=K_1\times_M K_2\times_M\ldots\times K_k/M$ 
and $h_{\bar w}[(k_1,\ldots,k_k)]=k_1\cdots k_ks_k\cdots s_1K_p$.
Let $\gamma_{\bar w}$ be a fundamental $\mathbf Z_2$-cycle in $\Gamma_{\bar w}$.

\begin{thm}[Bott -- Samelson]\label{thm:B-S}
$\{(h_{\bar w})_*[\gamma_{\bar w}]:\bar w\in W/W_p\}$ is a basis 
for the homology $H_*(K/K_p)$.  
\end{thm}

The relation between Theorems~\ref{thm:B-K} and~\ref{thm:B-S}
is given in the following theorem. The particular case 
when $G$ is a complex Lie group and the $K$-orbit 
is a principal orbit was already treated by Hansen in~\cite{Hansen}.  

\begin{thm}
The subset $W_u^{-1}\subset W$ contains precisely one representative
from each $W_p$-coset in $W$. For each $w\in W$, let 
$s_1,\ldots,s_k$ be as above and define $w_u=s_1\cdots s_k$
so that $w$ and $w_u^{-1}$ are $W_p$-conjugate. Then $w_u\in W_u$ 
and the image $h_{\bar w}(\Gamma_{\bar w})=K_1\ldots K_ks_k\ldots s_1
K_p$ of the Bott-Samelson cycle $h_{\bar w}:\Gamma_{\bar w}\to K/K_p$
is the closure of the Bruhat cell $N^+w_u^{-1}P_\Theta$. 
\end{thm}

\Pf A word about notation: 
if $\lambda\in\Sigma_r$ is a reduced restricted root, 
we write $\bar{\Lg}_{\lambda}=\Lg_{\lambda}+\Lg_{2\lambda}$.
(Of course it could happen that $2\lambda$ is not a 
restricted root in which case $\Lg_{2\lambda}$ would not be present.)

To see that $w_u\in W_u$ we need to verify that 
$\Sigma_{w_u}^+\cap\Span{\Theta}=\varnothing$. Since
the reduced restricted roots in 
$w_u^{-1}\Sigma^+\cap\Sigma^-$ are precisely 
$-\lambda_1,\ldots,-\lambda_k$,
that statement is equivalent to having 
$\{-\lambda_1^{w_u},\ldots,-\lambda_k^{w_u}\}\cap\Span{\Theta}
=\varnothing$, which follows from the fact that
$wp$ does not belong to $\ker\lambda_i$ 
for $i=1,\ldots,k$. Thus each $W_p$-coset 
in $W$ contains at least one element of $W_u^{-1}$. 
We invoke Proposition~1.1.2.13 in~\cite{Warner}
to see that the number of $W_p$-cosets in $W$ is equal to 
the cardinality of $W_u^{-1}$. This implies that 
each $W_p$-coset contains precisely one element of $W_u^{-1}$. 

Note that $U^+_{w_u^{-1}}$ is a nilpotent Lie group with Lie algebra
$\Lu^+_{w_u^{-1}}=\sum_{i=1}^k\bar{\Lg}_{\lambda_i}$ since
the restricted roots in 
$\Sigma^+_{w_u^{-1}}$ are $\lambda_1,\ldots,\lambda_k$. 
Let $U_i$, $i=1,\ldots,k$, denote the analytic subgroup of $G$ 
with Lie algebra $\bar{\Lg}_{\lambda_i}$.
Then $U_i$ is closed (because it is the connected component 
of the $N^+$-centralizer of 
the Lie subalgebra $\ker\lambda_i$ of $\La$)
and $U^+_{w_u^{-1}}=U_1\cdots U_k$ (because $U^+_{w_u^{-1}}$ is nilpotent), 
so Theorem~\ref{thm:B-K} implies that 
$N^+w_u^{-1}P_\Theta=U_1\ldots U_k w_u^{-1}P_\Theta$
and that the dimension of $N^+w_u^{-1}P_\Theta$ is equal to 
$\sum_{i=1}^k \dim(\bar{\Lg}_{\lambda_i})$. 

For each $j=1,\ldots,k$ define $w_j=s_j\cdots s_2s_1\in W$.
Set $B_j=w_jBw_j^{-1}$. Then $B_j$ is a Borel subgroup of $G$
with Lie algebra $\Lb_j=\Lm+\La+\sum_{i=1}^j\bar{\Lg}_{-\lambda_i}
+\sum_{i=j+1}^k\bar{\Lg}_{\lambda_i}$. Let also $P_j$  be 
the parabolic subgroup of $G$ with Lie algebra $\Lb_j+\bar{\Lg}_{\lambda_j}$.
Then we have the inclusions between groups
\[ P_1\supset B_1\subset P_2\supset B_2\subset P_3\supset
\ldots\subset P_k\supset B_k, \]
and we can form the quotient manifold
$\mathcal M_{w_u}=P_1\times_{B_1}\times P_2\times_{B_2}\cdots\times P_k/B_k$.

Since $K_j\subset P_j$ (because $P_j=K_jAw_jN^+w_j^{-1}$)
and $M\subset B_j$ (because $B_j=MAw_jN^+w_j^{-1}$), there is a
map induced by inclusion
$i:\Gamma_{\bar w}\to \mathcal M_{w_u}$ which is immediately seen to be 
injective and regular at $[(1,\ldots,1)]$ 
(since $K_j\cap B_j=M$). 
By equivariance $i$ is regular everywhere and it follows from the compactness
of $\Gamma_{\bar w}$ and the fact that 
$\Gamma_{\bar w}$ and $\mathcal M_{w_u}$ 
have the same
dimension that $i$ is a diffeomorphism. Hence we have the following 
commutative diagram
\[ \begin{array}{ccc}
     \Gamma_{\bar w} & \overset{\approx}{\longrightarrow}&\mathcal M_{w_u}\\
     \!\!\!\!\!\!\!h_{\bar w}\;\downarrow & 
     &\!\!\!\!\!\!\!f_{w_u}\;\downarrow \\
     K/K_p & \overset{\approx}{\longrightarrow} & G/P_\Theta \end{array} \]
where $f_{w_u}[(p_1,\ldots,p_k)]=p_1\ldots p_k s_k\ldots s_1P_\Theta$. 
Now $h_{\bar w}(\Gamma_{\bar w})=
f_{w_u}(\mathcal M_{w_u})=P_1\ldots P_ks_k\ldots s_1P_\Theta$
is compact and, as $P_j\supset U_j$,
 contains $U_1\ldots U_ks_k\ldots s_1P_\Theta=N^+w_u^{-1}P_\Theta$.
It follows that $h_{\bar w}(\Gamma_{\bar w})$ 
contains also the closure of $N^+w_u^{-1}P_\Theta$. 

We shall show next that we can write $P_j=U_jB_j\cup s_jB_j$ (disjoint union) 
where $U_jB_j$ is an open dense submanifold and $s_jB_j$ 
is a closed lower dimensional submanifold. 
For this purpose set~$U_{-j}=\theta(U_j)$ and
note that the $G$-centralizer $G_j$ of $\ker\lambda_j$
is a reductive Lie group, $\ker\lambda_j$ being a 
$\theta$-stable Abelian subalgebra of $\Lg$.
By the version of Theorem~\ref{thm:B-HC} for reductive groups
we can write $G_j=U_{-j}s_jMAU_{-j}\cup U_{-j}MAU_{-j}=
U_{-j}s_jMAU_{-j}\cup MAU_{-j}$, since $M$ and $A$ normalize $U_{-j}$. 
Now multiply through by $s_j$; we have $s_jU_{-j}s_j=U_j$
and so it follows that $G_j=U_jMAU_{-j}\cup s_jMAU_{-j}$. 
As we know, $K_j\subset G_j$ and $K_j$ acts
transitively on $P_j/B_j$ (with isotropy $M$). A fortiori,
$G_j$ acts transitively on $P_j/B_j$. Hence,
$P_j=G_jB_j=U_jB_j\cup s_jB_j$ as claimed. 

Owing to the above decomposition of $P_j$, 
by means of a straightforward induction argument 
one can easily establish that:
\begin{itemize}
\item[(i)] Any element in $\mathcal M_{w_u}$ can be represented by 
an element $(v_1,\ldots,v_k)\in P_1\times\ldots\times P_k$ where
$v_i\in U_i\cup\{s_i\}$. 
\item[(ii)] The image of $U_1\times\ldots\times U_k$ 
in~$\mathcal M_{w_u}$ is open and dense. 
\end{itemize}
Since $f_{w_u}$ is a diffeomorphism from $[U_1\times\ldots\times U_k]$
onto the Bruhat cell $N^+w_u^{-1}P_\Theta$, it finally follows 
from~(i) and~(ii) that $h_{\bar w}(\Gamma_{\bar w})$ equals 
the closure of $N^+w_u^{-1}P_\Theta$. \EPf

\bibliographystyle{amsplain}
\bibliography{ref}

\bigskip

\parbox[t]{7cm}{\footnotesize\sc Instituto de Matem\'atica e Estat\'\i stica\\
                Universidade de S\~ao Paulo\\
                Rua do Mat\~ao, 1010\\
                S\~ao Paulo, SP 05508-900\\
                Brazil\\ \hfill\\
                E-mail: {\tt gorodski@ime.usp.br}}{}\hfill{}
\parbox[t]{7cm}{\footnotesize\sc Mathematisches Institut\\
                Universit\"at zu K\"oln\\
                Weyertal 86-90\\
                50931 K\"oln\\
                Germany\\ \hfill\\
                E-mail: {\tt gthorber@mi.uni-koeln.de}}\hfil{}

\end{document}